\pdfoutput=1 
\documentclass[12pt]{article}
\usepackage{amsmath,amsfonts,amssymb,url,color,epsfig,graphics}

\newcommand{\R}{{\mathbb{R}}}

\newcommand{\Z}{{\mathbb{Z}}}

\newcommand{\Hy}{{\mathbb{H}}}
\newcommand{\To}{{\mathbb{T}}}

\def\ha{\frac{1}{2}}

\def\pa{\partial}
\def\ra{\rightarrow}
\def\ga{\alpha}

\def\gd{\delta}
\def\ge{\varepsilon}

\def\gg{\gamma}

\def\gl{\lambda}
\def\go{\omega}

\def\gs{\sigma}

\def\nghbd{neighbourhood~}

\newtheorem{lemm}{Lemma}[section]
\newtheorem{prop}{Proposition}[section]

\newtheorem{theo}{Theorem}[section]

\begin{document}

\title{Periodic geodesics for  contact sub-Riemannian 3D manifolds}

\author{Yves  Colin de Verdi\`ere\footnote{Universit\'e Grenoble-Alpes,
Institut Fourier,
 Unit{\'e} mixte
 de recherche CNRS-UGA 5582,
 BP 74, 38402-Saint Martin d'H\`eres Cedex (France);
{\color{blue} {\tt yves.colin-de-verdiere@univ-grenoble-alpes.fr}}}}
\maketitle
The goal  of this paper is to study periodic geodesics for sub-Riemannain metrics on a contact 3D-manifold.
We develop two rather independent subjects:
\begin{enumerate}
 \item  The existence of closed geodesics spiraling around periodic Reeb orbits for a generic metric.
 \item The precise study of the periodic geodesics for a right invariant metric on a quotient $\Gamma \backslash {\rm PSL}_2(\R)$ for which the Reeb flow is the geodesic flow of the
   corresponding hyperbolic surface $\Gamma \backslash \Hy $. 
   \end{enumerate}
In the first part (Section \ref{sec:setup} to Section \ref{sec:pb}), we   prove the following result  which was conjectured in  \cite{C-H-T-21}: 
\begin{theo} \label{theo:main} 
  Let $(M,D,g)$ be a contact 3D sub-Riemannian manifold
  and $\Gamma $ a periodic orbit of the canonically associated Reeb flow with period $T_0>0$. Let us assume that $\Gamma $ is non degenerate, meaning that $1$ is not
  an eigenvalue of the linearized Poincaré map of this orbit.  Then,  there exists a sequence 
  $\gg_k, ~k\geq k_0 $, of periodic sub-Riemannian geodesics of $(M,D,g)$ with 
  $\lim _{k\ra +\infty }\gg_k = \Gamma $ as closed sets with the Hausdorff topology and the length $l_k $ of $\gg_k$ is equivalent as $k\ra + \infty $ to $2\sqrt{\pi k T_0}$; more precisely,
  $l_k$ admits a full asymptotic expansion in powers of $k^{-\ha}$. 
\end{theo}
In fact,  as we will see, the $\gamma_k$'s are spiraling around $\Gamma $.
 The assumptions are generic for closed manifolds:  it is known that there exists closed Reeb orbits in any closed 3-manifold \cite{Ta-07}  and that they are generically non degenerate \cite{Bo-03}.
 The precise definitions of the terms in  Theorem \ref{theo:main}  will be given in Section \ref{sec:setup}.

 In the second part (Section \ref{sec:sl2r} to Section \ref{sec:omega}), we  give a description of the closed geodesics in the case of the {\it Liouville }
 contact structure on the quotient $\Gamma \backslash {\rm PSL}_2 (\R) $, the unit co-tangent bundle
 of an hyperbolic Riemann surface. Roughly speaking we show the existence of continuous families of 2D-tori on which the geodesics spiral linearly and are periodic for a dense set of parameters.
 This involves the Casimir Hamiltonian and the Euler equations in the dual of the Lie algebra of ${\rm PSL}_2 (\R) $.

 {\it Aknowledgements.--  Many thanks to Cyril Letrouit and Françoise Truc for their remarks allowing to improve the manuscript.
 Thanks also to Louis Funar and Christine Lescop for answering my questions about some  topological aspects. }

\section{Motivations}
The study of periodic geodesics is a classical part of  Riemannian geometry (see \cite{Be-03}, chap. 10).
There are not many works on closed geodesics on sub-Riemannian  manifolds, see however \cite{K-V-19,D-K-P-S-V-18} (for explicit calculations on spheres and Lie groups)
and \cite{Sh-21} and references therein  for the Heisenberg Kepler  problem. 
Starting from the work \cite{C-H-T-21}, it is natural to ask about closed geodesics spiraling around closed orbits of the Reeb flow.
A related motivation comes from inverse spectral problems: roughly speaking, is the set of periods of the Reeb flow  a spectral invariant of the sub-Riemannian Laplacian?
Our main theorem shows that it could  be true; indeed in \cite{Me-84}, Richard Melrose proved an extension of the
Chazarain-Duistermaat-Guillemin wave trace formula \cite{D-G-75} (see also \cite{CdV-07}) for sR contact Laplacian showing that the set of lengths of closed geodesics
(called the ``lengths spectrum'') is, generically,
a spectral invariant. One can then hope to recover the Reeb periods from the lengths of closed geodesics or at least to prove a rigidity result. 

\section{The setup} \label{sec:setup} 
A nice introduction to sub-Riemannian geometry can be found in the book \cite{Mo-02}, see also \cite{C-H-T-18} for what follows. 
Let us recall a few things: $M$ is a smooth manifold of dimension $3$, $D\subset TM$ is a smooth distribution of dimension $2$ defined globally by $D=\ker \ga $ where $\ga $ is a non vanishing real 1-form,
so that  $\ga \wedge d\ga $ is a non vanishing volume form. It follows that $M$ is oriented and $D$ is orientable. We choose some orientation of $D$.
The metric $g$ is a smooth metric  defined on $D$. We define the so-called co-metric $g^\star :T^\star M \ra \R^+$ by
\[ g^\star (q,p)=\| p _{|D_q} \| _{g(q)}^2 \]
where the norm is the norm on the dual of the Euclidean space $D_q$. 
The geodesic flow is the Hamiltonian flow of $\ha g^\star $ restricted to $g^\star =1$. The projections of the integral curves of the geodesic flow onto $M$ are the geodesics of
the sub-Riemannian manifold with speed $1$. 

To $g$ and the orientation of $D$ is associated a choice of a 1- form defining $D$ as follows:  we define $\ga _g $ so that
$\ker \ga_g =D$ and $d\ga _g $ restricts to $D$ as the oriented volume defined by $g$.

To $\ga _g $ is associated the Reeb vector field $\vec{R}$ on $M$
characterized by
\[ \ga_g (\vec{R})=1, ~d\ga _g (\vec{R},. )=0 \]
A periodic Reeb orbit is said to be {\it non degenerate}  if the linearized Poincaré map does not admit $1$ as an eigenvalue.

The Reeb vector field has the following Hamiltonian interpretation: the cone  $\Sigma =\{ (q,p)\in  T^\star M|p _{|D_q}=  0\}$,  generated by $\ga_g $,  is a symplectic sub-manifold of $T^\star M$.
We define the function $\rho (\ga ) =\ga/\ga _g $ on $\Sigma $. The function $\rho $ is homogeneous of degree $1$ and hence the Hamiltonian vector field $\vec{\rho}$ of $\rho $
is homogeneous of degree $0$. Let us denote $\Sigma^+:=\Sigma \cap \{ \rho >0 \}$ and by $\pi _\Sigma$ the projection of
$\Sigma^+ $ onto $M$. The projection  $(\pi_\Sigma)_\star(\vec{\rho})$ on $M$ is  well defined and is the Reeb vector field $\vec{R}$ (\cite{C-H-T-18}, sec. 2.4). 

Let us denote by $\pi $ the canonical projection of $T^\star M $ onto $M$. 
If $\Gamma $ is a periodic Reeb orbit, there exist a \nghbd $\Omega $ of $\Gamma $ and a 
 conical \nghbd $U$  of $ \pi_\Sigma ^{- 1}(\Omega  )$ in $T^\star M \setminus 0$   so that a  ``Birkhoff normal form'' holds in $U$.
This Birkhoff normal form is  defined as follows (see Section 5.1.4 of \cite{C-H-T-18}):
let us consider the conic symplectic manifold $\Sigma _\gs \times \R^2_{u,v}$ with the symplectic form \[ \go_\Sigma + dv\wedge du \] and the positive dilations
$\gl(\gs, u,v)=(\gl \gs, \sqrt{\gl}u,\sqrt{\gl}v) $. There exists an homogeneous symplectic diffeomorphism $\chi $ of $U$ onto an open  cone $V \subset \Sigma ^+\times \R^2$
so that,
$\forall \gs \in \pi_\Sigma ^{- 1}(\Omega  ), \chi (\gs )=(\gs, 0)$ and
\[ F(\gs, u,v):=g^\star \circ \chi^{-1}(\gs, u,v)=\rho I + \rho_2 I^2 + \cdots + O\left( I^2 (I/\rho)^\infty \right) \]
where $I=u^2+v^2$,  the $\rho_j$'s are functions homogeneous of degree $2-j$ on $\pi_\Sigma ^{- 1}(\Omega  )$ and
the remainder depends in general of $u$ and $v$ and not only of $I$ and $\gs$. Note that the remainder is natural: it is the expected estimate for a remainder which is flat along $\Sigma $ and homogeneous of degree
$2$.

In what follows, we will always study the geodesic flow in the Birkhoff coordinates. Note  that there exists some $I_0>0$ so that  the energy  shell $\{ F=1,~I< I_0 \} $ is
properly included in
the cone $V $ of $\Sigma \times \R^2$.  Hence, for $I $ small enough and $q \in \Omega $, we stay in the domain of the
Birkhoff normal form.

\section{The ``integrable'' case}\label{sec:int}

In this section, we will assume that the   Birkhoff normal form  is convergent: 
it means that $g^\star $  is symplectically equivalent in the cone $U$  to  some smooth function 
$ (\gs, I) \ra F(\gs, I )$ 
where  $F$ is a smooth homogeneous function  of degree $2$,  defined in the cone $V$.
Moreover, $F$ has  an asymptotic expansion as before \[ F(.,I) = \rho I + \rho_2 I^2 + \cdots \]
but the remainder depends only of $\gs $ and $I$. 
Clearly, in this case, the function $I$ is a first integral of the flow. Integrable here does not imply Liouville integrability  in general, because we have only two  integrals of the flow. 

Note that we will consider the geodesic flow (identified to the flow of $\ha F $) in the energy shell $\{ F=1, ~I<I_0 \} $.
For an Hamiltonian $H$ on a symplectic manifold, we denote by $\vec{H}$ the associated Hamiltonian vector field. 
\subsection{Closed geodesics} 
\begin{theo} Let us assume that there exists a periodic non degenerate Reeb orbit of period $T_0>0$ and that the Birkhoff normal form is convergent.
  Then, there exists  a sequence  of periodic  geodesics $\gg_k, ~k\geq k_0 ,$ of the sub-Riemannian manifold $(M,D,g)$ accumulating on $\Gamma $ with
  lengths 
  \begin{equation} \label {equ:lengths}
    l_k= 2\sqrt{\pi k T_0} + \sum _{j=0}^\infty a_j k^{-j/2} + O\left(k^{-\infty }\right) \end{equation}
\end{theo}
{\it Proof.--} Let us denote by $H_I$ the Hamiltonian on $\Sigma $ defined by $\ha F(., I)=\ha \rho I + O(I^2)$. For the Hamiltonian $H =\ha g^\star $ of the geodesic flow expressed in  the ``Birkhoff coordinates'', we have
\[ \ha \vec{F}=\overrightarrow{H_I} + \frac{\pa F}{\pa I}\pa _\theta \]
where
\begin{equation}\label{equ:polar}
  (u,v)=(\sqrt{I}\cos \theta , \sqrt{I}\sin \theta ).
\end{equation}

Let us start with the
\begin{lemm} Let us consider the map $\pi_I $ which is the restriction of $\pi_\Sigma $ to $\{ \gs| F(\gs , I)=1 \}$.
  Then, for $I$ small enough, $\pi_I$ is a  diffeomorphism over a fixed \nghbd of $\Gamma $.
\end{lemm}

{\it Proof.--}
For $I$ small enough, there exists a smooth function $\gl _I : \{ \rho =1\} \ra \R^+ $
so that $F(\gl_I (\gs) \sigma, I)=1$. The function $\gl_I $  admits an expansion $\gl_I= 1/I + O(1)$.
We define $\Lambda _I:\{ \rho =1\} \ra \{ F(., I)=1\} $ by $\Lambda _I (\gs)=\gl_I (\gs ) \gs $.
The map $\Lambda _I $ is a diffeomorphism. We can hence consider the map $\pi_I \circ \Lambda _I $ from $\{ \rho =1\} $ on $M$.
We have $\pi_I \circ \Lambda _I (\gs )= \pi_I (\gl_I (\gs) \gs )= \pi _\Sigma (\gs )$ which is clearly a diffeomorphism. \hfill $ \square $

We consider then the geodesic flow. For each value of $I$, we project it on $\Sigma $, i.e. we consider the flow of the Hamiltonian $\ha F(.,I)$ on $\Sigma $ restricted to $F=1$.
We have the
\begin{lemm} The previous flow, projected by $\pi_I$ and with a change of time $s=2t/I $, is a smooth perturbation of the Reeb flow on $M$.
\end{lemm}{\it Proof.--}
The flow projected on $\Sigma $ is $\ha I \vec{\rho  } + O(I^2)$. Hence the change of times reduces to $\vec{\rho } + O(I)$.
This last vector field projects on $M$ as $\vec{R}+O(I)$. \hfill $ \square $

From this and the fact that $\Gamma $ is non degenerate (see Appendix \ref{app:poincare}), we obtain a periodic  orbit of the Hamiltonian $F(.,I)$ of period $T(I)$ with

\begin{equation} \label{equ:T(I)}  T(I)\sim \frac{2T_0}{I} +\sum_{j\geq 0} c_j I^j  \end{equation} 

We have now to close the angular part of the dynamics given by $\pa F/\pa I \pa_\theta $.
\[ \theta (T(I) )-\theta (0)= \int _0^{T(I)} \frac{\pa F}{\pa I}(\Gamma _I(t),I ) dt  \]
The righthandside  of this  equation admits a full expansion
\[ \frac{2T_0}{I^2}\left(1+\sum_{j=1}^\infty  b_j I^j  \right) \]
which has to be equal to $2k \pi $ in order to close the  geodesic.
Hence, we get an  asymptotic expansion of $I$ in terms of  powers of $k^{-\ha}$ which gives the asymptotic of the lengths by inserting into the Equation (\ref{equ:T(I)}).
\hfill $\square $ 

\subsection{Poincaré section}
Let us now describe a  {\it Poincaré section} $S_0$ and the corresponding {\it Poincaré map} $P_0$ for the geodesic flow  assuming for simplicity that $F=\rho I$.
For the definitions  and properties of the Poincaré maps, see  Appendix \ref{app:poincare}. 

Let  $X_R $ be a Poincaré section of the Reeb flow in the energy shell  $\rho =R$ so that
all $X_R$'s project on a fixed Poincaré section of the periodic orbit of the Reeb flow in $M$. Let us  use polar coordinates given by Equation (\ref{equ:polar})
identifying  $\R^2 _{u,v}\setminus 0$ to $T^\star _{\theta, I}(\R/2\pi \Z)$.
The latter manifold will be denoted by $T^\star $ in what follows. 
The manifold 
\begin{equation} \label{equ:Poincare}  S_0:= \{(\gs,\theta, I)| \gs \in  X_R, ~ RI=1, ~ I< I_0\}  \end{equation} 
  is  a Poincaré section of the geodesic flow.
 But the Reeb flow projects onto $M$; we can hence 
identify any   Poincaré section $X_R$   with a Poincaré section  of the Reeb flow in $M$ denoted by $X$, which can be assumed to be independent of $R$. 
    The   Poincaré $S_0$ section is then parametrized  by ${\bf S}_0$ given by
    ${\bf S}_0=   X_q  \times T^\star_{\theta, I} $. 
    The Poincaré map is  given in these coordinates by
    \[ P_0(q, \theta ,I )= (\Pi (q), \theta + 2T(q)/I^2, I ) \]
    where $T(q)$ is the return time of the Reeb flow when starting from $q$ and $\Pi$ is the Poincaré map of the Reeb orbit. 
    Note that the  cylinder $C_0:= \{  q_0 \} \times T^\star$ is invariant by $P_0$ as well as the symplectic form $d\theta \wedge dI$ on it.  Note also that the latter symplectic form is the
    restriction of the symplectic form on $\Sigma \times T^\star$ to $C_0$. 
    We will show latter  that the Poincaré map is weakly perturbed when  $I$ is small in the non integrable case. 

\section{Summary of the proof}

In what follows,  we will assume, for technical simplicity,  that the co-metric admits the simple normal form $g^\star = \rho I + O\left(I^2 (I/\rho)^\infty \right) $
in some conical \nghbd of $\pi^{-1}(\Gamma) \cap  \Sigma ^+$. 
The general case would be $g^\star =F(\gs, I) + O\left(I^2 (I/\rho)^\infty \right) $ with $F$ a sum of the BNF as in the previous section. The same strategy will work in the latter  case.
So we will see the flow of $\ha g^\star $ as a perturbation of the flow of $\ha \rho I  $. 

The scheme of the proof is as follows.
\begin{enumerate}

\item Using  the structural stability of the non degenerate closed Reeb orbit, we get  a symplectic cylinder $C $ of dimension $2$ close to $C_0$ invariant by the Poincaré map.

\item We show the existence of invariant circles $c_k $ inside that cylinder. 

\item We  apply the Poincaré-Birkhoff fixed point theorem  to the annuli between $c_{k-1}$ and $c_{k+1}$. 


\end{enumerate}

\section{The invariant cylinder $C$}

Let us first build a Poincaré section for the geodesic flow. Let
$S_0\subset \Sigma \times T^\star $ be the Poincaré section for $\ha \rho I$ given in Equation (\ref{equ:Poincare}). 
We choose a Poincaré section which is a dilation of $S_0$ (respecting the fibers) of the form
$\gl  S_0$ with $\gl :\Sigma \times T^\star \ra \R^+$ a germ along $\pi^{-1}(\Gamma)$ . We have to satisfy
$g^\star =1$ with $g^\star =\rho I + O((\gl I)^\infty )$ and $\rho I =1$ on $S_0$ , this gives
\[ \gl ^2+ O((\gl I)^\infty ) =1\]
and hence $\gl =1+O(I^\infty )$. This is clearly a Poincaré section because it is still transversal to the flow.
Using the projection onto $M$, we can still identify this 
Poincaré section with  ${\bf S}_0$.
We  need the
\begin{lemm}\label{lemm:geod}
  On the energy shell $F =1$ and for times $O(1/I)$ the geodesic flow differs from the unperturbed one by $O(I(0)^\infty )$. 
\end{lemm}
{\it Proof.--}
We put $I_0:=I(0)$.
First, from $dI/dt = O(I^\infty )$, we get $I(t)=I_0+O (I_0^\infty )$ for times $t=O(1/I)$. We get  $\rho =1/I_0+O (I_0^\infty )$.
The Lemma  follows then by looking at
\[\ha  \overrightarrow{ g^\star }= \frac{I}{2}\vec{\rho }+ \rho \pa_\theta + O (I^\infty )\]
where we can replace $I$ by $I_0$ and $\rho $ by $1/I_0$ modulo $O (I_0^\infty )$.
We use then Equation (\ref{equ:perturb}) in Appendix \ref{app:perturb}.
Recall that the flow of $\ha \vec{\rho I}$ is given by
$G_t^0 (\gs, \theta ,I )=\left( R_{It/2}(\gs), \theta + \frac{t}{I}, I  \right) $.
The differential of that flow is hence
\[ DG_t^0 (\gs, \theta ,I )=
  \begin{pmatrix} DR_{It/2} & 0 &0 \\
    0 &1&0 \\
    0 &-\frac{t}{I^2}&1 
  \end{pmatrix} \] 
We get the following estimate
\[ \|\left( DG_t^0 (\gs, \theta ,I )\right)^{-1}\| \leq \| \left( DR_{It/2}(\gs)\right)^{-1} \|+ |t/ I^2 | \]
Hence, for times $O(1/I)$, we have
\[ \|\left( DG_t^0 (\gs, \theta ,I )\right)^{-1}\| =O\left(I^{-3} \right) \] 
Using Equation (\ref{equ:perturb}) in Appendix \ref{app:perturb} and the notations there, we get
$ \|\frac{d}{dt}"w(t,x)"\|= O\left(I_0^\infty \right)$ for times $O(1/I_0)$ and $ w(t,x)=x+ O\left(I_0^\infty \right)$. The result follows.
\hfill $\square $

We get
\begin{prop} For $a>0$, let us denote $T^\star _a :=T^\star \cap \{ |I|<a \}$.
  There exists $a>0$ and smooth functions $q:T^\star _a \ra X $ and $T:T^\star _a \ra \R $ so that
  \begin{enumerate}
  \item \[ q(\theta, I) =q_0+O\left(I^\infty \right)\]
  \item \[ T(\theta, I)= \frac{2T_0}{I} + O\left(I^\infty \right)\]
    \item 
  The flow $G_t$ of $\ha g^\star $ restricted to $g^\star =1$ satisfies
  \[ G_{T(\theta , I )}( q(\theta , I ), \theta ,I)= \left(q(\theta , I ), \theta + \frac{2T_0}{I^2}+O\left(I^\infty \right) , I+ O\left(I^\infty  \right) \right)\]
\end{enumerate}
\end{prop}
Hence the cylinder \[ C:=\{ (q(\theta, I), \theta ,I) |(\theta ,I)\in T^\star \} \]
is invariant by  $ G_{T(\theta , I )}$ which is the restriction of the Poincaré map  to $C$. Moreover the symplectic form restricts to $C$ as a symplectic form $dI \wedge d\theta + O\left(I^\infty \right)$.

{\it Proof.--}
Let us consider the return map to $X$ of the $q$-component for $(\theta, I)$  fixed with $I$ small. It follows  from Lemma \ref{lemm:geod} that the return map is
$O\left(I^\infty  \right)$ close to the unperturbed Poincaré map of the Reeb flow with a 
return time  $O(1/I)$. The conclusion  follows then from  the non degeneracy of the Reeb orbit.
\hfill $\square $

    \section{The invariant circles $c_k$}

    Applying Theorem \ref{theo:inv}  to the restriction $P_C$ of $P$ to the cylinder $C$ and the circles $c_k^0:= \{ (\theta ,I)| I = \sqrt{T_0/k\pi } \}$, we will get circles $c_k $ globally, but not
    pointwise, invariant by the map $P_C$
    and close to $c_k^0$.

    More precisely,
     the restriction $P_C$ of the Poincaré map to the cylinder $C$ writes
     \[ P_C(\theta, I )= (\theta +2T_0/I^2, I )+ O\left(I^\infty \right)\]
     using the coordinates $(\theta, I)$ in $C$. 
    Putting $J=1/I^2$ and $\theta '= \theta +2k\pi  $ , we get the  map:
    \[ S(\theta', J)=(\theta '+ (2T_0 J-2k\pi), J ) +O\left( J^{-\infty } \right) \]
    Near $J=k\pi /T_0$, we  get
    \begin{equation}\label{equ:S}
      S(\theta', J)=(\theta' + (2T_0 J-2k\pi), J ) +O\left(k^{-\infty} \right)\end{equation}
    This map is a perturbation of the map
    $(\theta' ,J)\ra (\theta' + (2T_0 J-2k\pi), J )$ to which we can apply  Appendix \ref{app:inv}
    with $X=T^\star $ and $Y_0=\{ J=k\pi /T_0\} $ and get the curves $c_k $ which are  $O\left(I^\infty \right)$ close to $c_k^0$.

    \section{Applying the Poincaré-Birkhoff theorem}\label{sec:pb}
    
    We  apply the 
    Poincaré-Birkhoff Theorem \ref{theo:PB}.   to the annuli between $c_{k-1}$ and $c_{k+1} $ and the lift to the universal cover of the annulus 
    moving the lift $C_{k+1}$  of $c_{k+1}$ by a map close to $\theta \ra \theta +2\pi $
      and the lift $C_{k-1}$  of $c_{k-1}$ by a map close to $\theta \ra \theta -2\pi $.
      This way we get a fixed point in the annulus for the lift. This point is $O\left( I^\infty \right) $ close to $c_k^0$,  because the other points are moved at some speed like $I^m$ for some $m \geq 0$: it is because
      of the estimate  (\ref{equ:S}) of the map $S$.
This finishes the proof of Theorem \ref{theo:main}.

\section{Periodic geodesics on $\Gamma \backslash  {\rm PSL}_2(\R) $}\label{sec:sl2r}
The goal of this second part, which is quite independent of the previous ones, is to describe the periodic geodesics of a specific right invariant sub-Riemannian contact structure on compact quotients
\[\Gamma \backslash  {\rm PSL}_2(\R) :=\{ \Gamma.g |g\in {\rm PSL}_2(\R) \}\]
We assume for simplicity that $\Gamma $ has no elliptic elements, so that all elements of $\Gamma \setminus {\rm Id }$ are hyperbolic.  Because we took the quotient ${\rm PSL}_2(\R)$ 
of ${\rm SL}_2(\R) $
by $\pm {\rm Id}$, we can represent any hyperbolic transform by a matrix with eigenvalues $\gl, ~1/\gl $ with $ \gl >1$.

Our analysis  can be extended to some Riemannian case (see \cite{Sa-98}) or to another sub-Riemannian structure like the magnetic one (see \cite{Ch-20}), or even to any right invariant Hamiltonian.
Note also that the Quantum version of this study, namely the spectral theory of the associated sub-Riemannian Laplacian, follows some parallel path (see \cite{C-H-W-22?}). 
\section{Lie-Poisson bracket}

Let $G$ be a Lie group. We identify its Lie algebra ${\mathcal G}$ to  the space of right invariant vector fields equipped with the
bracket of vector fields. We consider also  the algebra of right invariant differential operators, called the envelopping algebra,
generated by the Lie algebra.
The principal symbols of these right invariant operators  are determined by their values on the dual ${\mathcal G}^\star =T^\star _{\rm Id}G$
of the Lie algebra. This gives a Poisson bracket on ${\mathcal G}^\star $.  We need only  to compute
it for coordinates functions: 
if $X\in {\mathcal G}$, the symbol of $X$ at $g={\rm Id}$ is $\gs (X) (p)= ip(X) $.
If $X,Y \in {\mathcal G}$, 
the symbol of $[X,Y]=Z$ is $-i \{ \gs (X), \gs (Y) \} $.
Hence the formula 
\[ \{ a, b \} (p )= -[ da (p),db (p )] (p)  \] 
can be checked for operators of the Lie algebra:
if $p\in {\cal G}^\star$, 
\[ \gs ([X,Y])(p)=-i  \{ \gs (X), \gs(Y )  \}(p)=  i p ( [ X,Y]|)   \]
or, if the functions $\xi $ and $\eta $ on ${\cal G}^\star$ are defined by $ \xi(p)=p(X), ~\eta (p)=p(Y)$, 
\[ \{ \xi, \eta \} (p)=-p([X,Y])\]
This bracket is called the {\it Lie-Poisson bracket}. The associated Hamiltonian dynamics are given by {\it Euler equations.}

This gives the  {\it Euler equations} using for  $a$ in the previous equation a  coordinates system on   ${\mathcal G}^\star $.
All of this is explained in the book \cite{M-R-98}, Sec. 13.1.

Let us compute the Poisson bracket in the case of $G={\rm PSL}_2(\R)$.
The Lie algebra is the 3D space of trace free real $2\times 2$ matrices
\[ M(x,y,z):=\left( \begin{matrix} z & x \\ y & -z \end{matrix} \right) \] 
We write $M= xX  + yY +zZ$. We have
\[            [X,Y]=-Z,~[X,Z]=2X, ~[Y,Z]=-2 Y \]
There is a minus sign w.r. to the matrix bracket because of the right invariance!

Hence, if $(\xi,\eta,\zeta )$ are the coordinates on ${\cal G}^\star$ dual to $(x,y,z)$, we have
\[ \{ \xi, \eta  \}=\zeta,~  \{ \xi, \zeta  \}=-2\xi,~    \{ \eta, \zeta  \}=2\eta \]
The symplectic leaves of the Poisson bracket are  the connected components of  the level sets of the
{\it Casimir Hamiltonian} defined by ${\rm Cas}(\xi,\eta,\zeta ):=\ha\zeta^2 +2 \xi \eta $.
The  Hamiltonian ${\rm Cas}$ Poisson commutes with all functions on ${\cal G}^\star $.

We consider the right invariant contact distribution generated by $X$ and $Y$ and the sub-Riemannian metric $g$ for which  $(X,Y)$ is an orthonormal basis.
The right invariant Hamiltonian of the geodesic flow is $\ha g^\star $ with  $g^\star =\xi^2+\eta ^2 $.

\section{Invariant 2-tori, the frequency $\go $ and the period $\tau (C)$}
The Hamiltonian $g^\star $ Poisson commutes with ${\rm Cas}$. Both flows are complete because the momentum map
$({\rm Cas},g^\star):T^\star M \ra \R^2$ is proper. We get hence an $\R^2$-action $\Phi $ on $M\times{\cal G}^\star \equiv T^\star M$ preserving the values of the momentum map. 
We will see that for most values of $C$, there exists  orbits of $\Phi $ which are tori $\To_C$ on which the action $\Phi $ is linear and such that the orbits of $\overrightarrow{{\rm Cas}}$
are periodic of period $T^{\rm Cas}(C)$.
The geodesic dynamics on
$\To_C $ induces a Poincaré map on an orbit of the Hamiltonian vector field $\overrightarrow{{\rm Cas}}$ of ${\rm Cas}$. This map is a translation of $t\ra t +\go (C) /T^{\rm Cas}(C) $ of the choosen Casimir orbit.
For $\go =p/q~ {\rm mod} \Z$ with $p,q,~q\geq 1,$ which are  coprime,  we get periodic geodesics. These geodesics are crossing  $q$  times  the Casimir periodic orbit.

For $C\in ]-1, +\infty [\setminus \{ 0,1 \} $, the geodesic Hamiltonian $\ha g^\star $ admits  two periodic orbits  on the symplectic leaf ${\rm Cas}^{-1}(C)$. These orbits are supported by
the two connected components of $\left( g^\star \right)^{-1}(1)\cap {\rm Cas}^{-1}(C)$. We denote by
$\tau (C)>0$ the period of these orbits.

Note also that the antipodal map $p\ra -p $ changes the orientation of both dynamics. 

\section{Main result} 
To each closed sub-Riemannian geodesic $c$ of period $T>0$, we associate the following invariants:
\begin{itemize}
\item The {\it free homotopy class}
  of the projection $\pi (c) $ onto $\Gamma \backslash \Hy $ which is given by a conjugacy class $[\gamma ]$ in $\Gamma $.
  We denote it by $F_\Hy (c) $. 
\item The {\it spiraling integer}  of any horizontal smooth closed curve in $M$ of period $T$  is defined as follows:
  the right invariant distribution $D$ spanned by $X,Y$ is a trivialized bundle. We look at the
  map $\theta :\R/T \Z \ra \{ \xi^2+\eta ^2 =1 \} $ which associates to $t$ the vector $(\xi, \eta )$ so that $\xi (t) X (c(t))+\eta  (t)Y(c(t)) $ is  the tangent vector to the curve  at the point of parameter $t$. 
  The spiraling integer is the degree of that mapping.
  Note that this integer is invariant by $C^1$ homotopies of horizontal periodic curves.  We denote it by ${\rm Sp}(c)$.
\item The {\it Morse index} $\mu (c)$ which is the Morse index of the energy functional restricted to horizontal curves.
\item The {\it free homotopy class} of $c$ as a conjugacy class in $\pi_1 (M)$ denoted by $F_M (c)$.
 \end{itemize}
We want to prove the following
\begin{theo}   For each value of $C_0\in [-1, +\infty [\setminus \{0 \} $ there exists closed sub-Riemannian geodesics whose momentum satisfies
  ${\rm Cas}(p)=C_0$, $g^\star =1$. More precisely: 
  
\begin{itemize}
\item If $C_0>1$,  for each primitive hyperbolic $\gg \in \Gamma $, there exists a 2D-torus $\To_\gg ^{C_o} $ invariant by the geodesic flow. The dynamics of the geodesic flow on 
  $\To_\gg^{C_o} $ is linear and periodic for the  dense set of values of $C_0$ for which  $\go (C_0)= p/q$ is rational. For such closed geodesics, we have 
 ${\rm Sp}(c)=q$ and the length $l(c)=q\tau (C_0)$. 
\item   If $C_0=1$,  for each primitive hyperbolic $\gg \in \Gamma $, there exists one   periodic geodesic of length $\sqrt{2} \log \gl  $ with $F_\Hy (c) =[\gg ]$
  and ${\rm Sp}(c)=0$. 
  \item If $0< C_0 <1$,  for each primitive hyperbolic $\gg \in \Gamma $, there exists a 2D-torus $\To_\gg ^{C_o} $ invariant by the geodesic flow. The dynamics of the geodesic flow on 
    $\To_\gg^{C_o} $ is linear and periodic for a dense set of values of $C_0$ with $\go(C_0)=p/q$. For such closed geodesics, we have 
   ${\rm Sp}(c)=0$ and the length $l(c)=q\tau (C_0)$. 
\item   If $-1< C_0<0$, there exists a 2D-torus $\To_{\rm Id}^{C_o}  $ invariant by the geodesic flow.  The dynamics of the geodesic flow on 
$\To_{\rm Id}^{C_0} $ is linear and periodic for a dense set of values of $C_0$ for which  $\go (C_0)= p/q$ is rational. We have   $F_\Hy (c) =\{ \rm Id \}$, ${\rm Sp}(c)=0$ and the length $l(c)=q\tau (C_0)$. 
\item   If $C_0=-1$, there exists one  periodic orbit of length $2\sqrt{2} \pi   $ with $F_\Hy (c) =\{ {\rm Id }\} $
  and ${\rm Sp}(c)=0$. 
\end{itemize}

\end{theo}
It is not clear for us how to compute the other invariants...

\section{Casimir periodic orbits}
Recall that  $\Gamma $ is a co-compact lattice  with no elliptic elements in $G$ and $M= \Gamma \backslash G$.
The compact 3-manifold $M$ can be identified with the unit cotangent bundle of the compact Riemann surface $\Gamma \backslash \Hy $ where $\Hy $ is the Poincaré half-plane
(see Appendix \ref{app:hyp}). 
We consider the Hamiltonian ${\rm Cas} $ on $T^\star M\equiv M \times {\cal G}^\star $. 
We have
 \[ \overrightarrow{{\rm Cas}} = \zeta \vec{\zeta}+2 \xi \vec{\eta } + 2 \eta \vec{\xi} \]
On the other hand
$(\xi,\eta,\zeta )$ are first integrals, hence
for $(\Gamma .g ,p)\in M\times {\cal G}^\star $ with $p=(\xi, \eta,\zeta)$, if ${\rm Cas}_t$ is the Casimir flow, we have
\[ {\rm Cas}_t ( \Gamma . g, p)= (\Gamma . ge^{tA(p)}, p) \] with
\[ A(p):= \zeta Z + 2\xi Y +2\eta X \]
and ${\rm Cas}=-\ha {\rm det }(A(p))$.

We get a periodic orbit $t\ra (\Gamma .ge^{tA(p)},p)$  of period $T$ if and only if 
$\Gamma .g e^{TA(p)}=\Gamma .g$, i.e. 
    \[ ge^{TA(p)} g^{-1} =\gamma \]
    for some $\gamma \in \Gamma $ and $g\in G$. That periodic orbit is the left translate of the invariant line of $\gamma $ by $g^{-1}$. 
    
\begin{itemize}
  \item 
   If $\gamma $ is hyperbolic with eigenvalues  $\gl$ and $ 1/\gl $ with 
  $\gl >1$, 
   we need to have
  that the eigenvalues of $A(p) $ are real and non zero, hence $C>0$ and
  then
  \[ e^{T \sqrt{2C}}=\gl \]
  or
  \[ T^{\rm Cas}(C)= \frac{1}{\sqrt{2C}}\log \gl \]
    Note that this orbit is simply a translate of a lift of a periodic geodesic of $\Gamma \backslash \Hy$. 
  
\item
  If $-1\leq C<0$, then $e^{TA(p)} $ is a rotation of angle $T \sqrt{-2C}$ and hence periodic of period $T^{\rm Cas}(C)=2\pi /\sqrt{-2C}$, independently of $\Gamma $.
  This orbit is homotopic to the the compact subgroup ${\rm SO}_2(\R)$ of $G$. 
  \item If $C=0$, we get  no periodic orbits. 
  \end{itemize}

\section{The geodesic flow and periodic geodesics}\label{sec:gf}

The sub-Riemannian geodesic flow $G_t$, which is the Hamiltonian flow  of  $\ha g^\star $, acts on the set of Casimir periodic orbits of fixed period:
if ${\rm Cas}_T (z)=z$, we have $ G_t({\rm Cas}_T(z))= {\rm Cas}_T( G_t(z) )$, hence
${\rm Cas}_T (G_t(z)) =G_t(z) $.
We consider  the sub-Riemannian geodesic flow   $G_t$ in the energy shell $g^\star =1$ giving geodesics with speed $1$. 
Recall that the  periodic orbits  of $\overrightarrow{{\rm Cas}}$  are parametrized by the value of $p(0)$. 
We can hence look  at the Poisson action of $\ha g^\star $ on ${\cal G}^\star $ which preserves the symplectic leaves.  
For each values $C_0$ of ${\rm Cas}$, this Poisson  action  is  an Hamiltonian action on the  2D symplectic manifold  ${\rm Cas}=C_0$ with Hamiltonian ${\ha (\xi^2+\eta^2)}$.
We parametrize the cylinder $\{  \xi^2+\eta^2=1 \} \subset {\cal G}^\star $ by
$(\xi =\cos \theta, ~\eta=\sin \theta,\zeta )$ with $\theta \in \R/2\pi \Z$ and $\zeta \in \R$. We get
$ {\rm Cas}= \ha \zeta^2+ \sin 2\theta  $.

Each time this Poisson  action on ${\cal G}^\star $ is periodic, the  orbit of the $\R^2$-action $\Phi$  is a torus $\To_C$ on which $\Phi $ acts linearly. 
To each such a torus, we associate a rotation  number $C\ra \go(C) \in\R/\Z$ as follows: a periodic orbit of the Casimir flow of period $T$
is a Poincaré section of the geodesic flow $G_t$ on that torus. The frequency $\go$ is the rotation number of the Poincaré map which is conjugated to a rotation.
If $\go =p/q $ is rational with , we get a 1-parameter family of closed geodesics. We will sometimes consider a lift $\tilde{\go } $ of $\go $ to $\R$. 
We will in particular study  below the map $\tilde{\go} :]-1, +\infty [\setminus \{ 0,1 \} \ra \R$.
The lengths of these periodic geodesics is $q\tau (C)$ because they have to spiral $q$ times before closing.

\subsection{Level sets of ${\rm Cas}$ on $\xi^2+\eta ^2 =1$}
The geodesic flow preserves the level sets  of ${\rm Cas}$ restricted to the unit bundle. 
We have hence to look at the different types of level lines. Due to right invariance, we have only to  look at the restriction of these sets to $Y:={\cal G}^\star \cap \{ g^\star =1 \}$ which is a 2D cylinder.
The critical values of ${\rm Cas}$ restricted to $Y$ are
the global minimum $-1$  and a saddle point $+1$. We have also to look at the special value ${\rm Cas}=0$ separating the elliptic from the hyperbolic case for the Casimir flow.  

\includegraphics[scale=0.5]{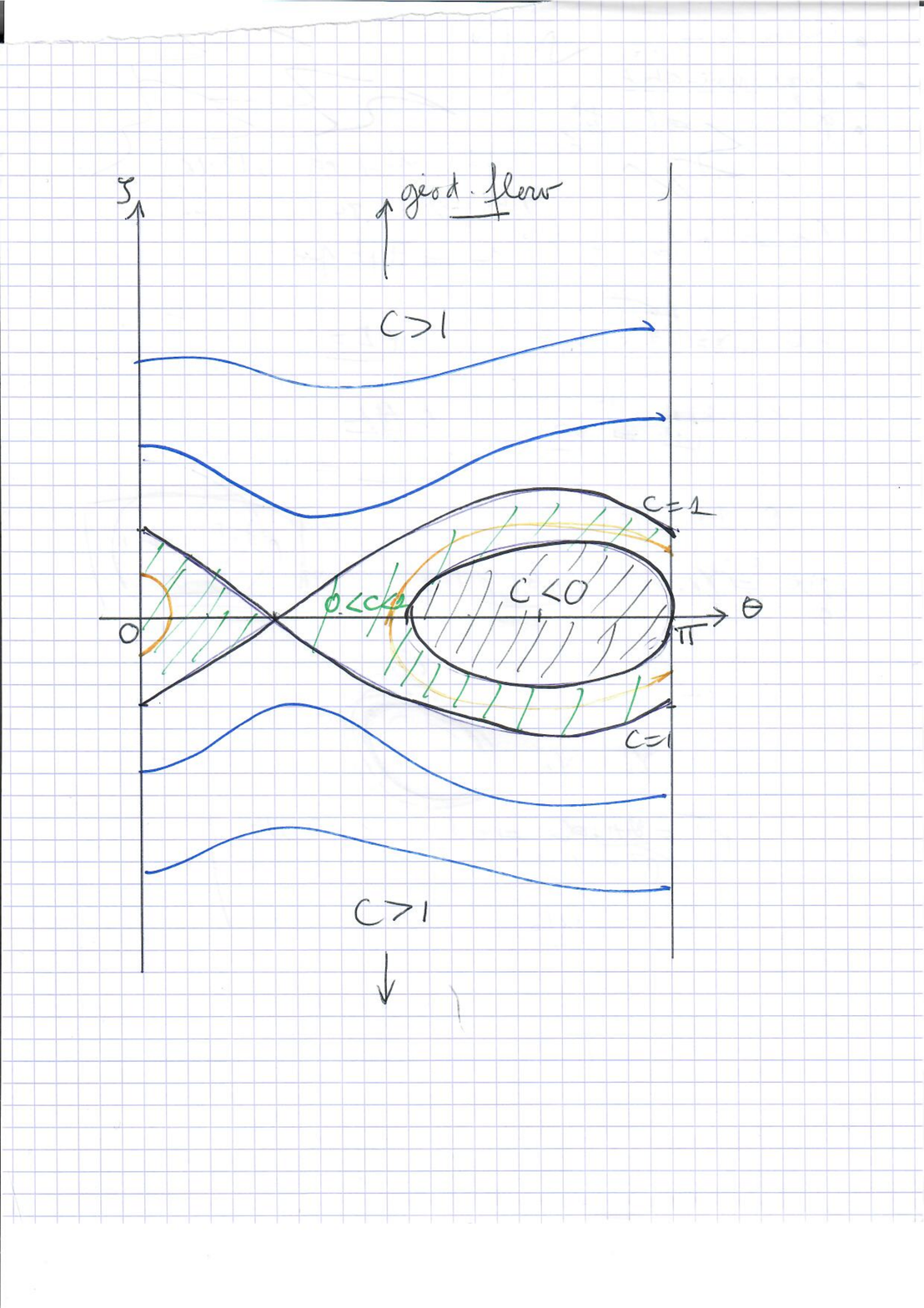}

We have hence many cases to consider which are the classical versions of the  irreducible representations of ${\rm PSL}_2 (\R )$ (see \cite{Tay-86}, chap. 8). 
\subsubsection{The case  $C_0>1$: ``principal series''}

The action of  $\ha g^\star $ on the Poisson manifold ${\cal G}^\star$   restricted   to $g^\star =1$ is  periodic.
This implies that the corresponding orbits of the $\R^2$-action is a 2-torus. We get hence periodic geodesics each times $\go $ is rational.

\subsubsection{The case where  $C_0=1$}

In this case, the orbits of $\Phi $  are cylinders and circles. The circles corresponds to the critical points of ${\rm Cas}$ with critical value $1$. 
For the critical point $\zeta =0, \xi=\eta =1/\sqrt{2}$, we have  periodic geodesics parametrized by the hyperbolic conjugacy classes of $\Gamma $
with lengths the lengths of the corresponding Riemannian periodic geodesics of $\Gamma \backslash \Hy $. 

\subsubsection{The case where  $0< C_0< 1$: ``complementary series''}
Again we get tori as orbits of the $\R^2$-action and periodic geodesics each times $\go $ is rational.
\subsubsection{The case where  $C_0=0$}
No periodic orbit  for the Casimir flow which is parabolic. 

\subsubsection{The case where  $-1< C_0< 0$: ``discrete series''}
Again we get tori as orbits of the $\R^2$-action and periodic geodesics each times $\go $ is rational.

\subsubsection{The case where  $C_0=-1$}For the critical point $\zeta =0, \xi=-\eta =1/\sqrt{2}$, we get  periodic geodesics which are the orbits of the
${\rm SO}_2 $ action, i.e. the fibers of $M\ra \Gamma \backslash \Hy $. 

\section{The frequency map $\go $}\label{sec:omega}

If $\go (C)$ is rational, we get a family of periodic geodesics. We will
 study the function $\tilde{\go} $ and show that it is rational for a dense set of values of $C$. 
We have to study first two  other functions of $C$ the period $T^{\rm Cas}$ of the Casimir flow and the
period $T^{\rm geod}$ of the action of the geodesic flow on the periodic orbits of the Casimir flow.
\subsection{The function $T^{\rm Cas}$ }
We know the that $T^{\rm Cas}(C) =l_\gg /\sqrt{2C}$ for $C>0$  and
$T^{\rm Cas}(C) = \sqrt{2}\pi / \sqrt{-C} $ for $C<0$. 
\subsection{The function $T^{\rm geod}$}
The function $T^{\rm geod}$ is smooth and non vanishing on $]-1,+\infty [ \setminus \{1 \}$.
For  $C\ra  +\infty $ is $T^{\rm geod}\sim 1/\sqrt{C}$, the limits at $+1$ are  $+\infty $, the limit at $-1$ is some $>0$ number. 
\subsection{The function $\go$}
Recall that the function $\tilde{\go} $ is a lift to $\R$ of the rotation number of the Poincaré map induced on the periodic Casimir orbit by the geodesic flow. If $\go (C)$ is rational, we get a full
torus of periodic geodesics. In order to understand $\go $, we look at the projection of the tori on $M$. The projected Hamiltonian vector fields are
respectively $V^{\rm Cas}= 2\eta X+2\xi Y +\zeta Z$
and $V^{\rm geod}= \xi X + \eta Y $. Both are independent outside the critical points $\zeta =0, \xi=\pm \eta = \pm 1/\sqrt{2}$. We  consider the angle
$\ga $ between the two vectors (in the Riemannian metric on $M$ whose $(X,Y,Z)$ is an orthonormal basis).
We have
\[ |\cos \alpha |= \frac{|\xi \eta |}{ \sqrt{1+\zeta ^2/4}} \]

Hence, we get the following behaviour of $\go $:
\begin{itemize}
\item As $C \ra +\infty $, both periods are of the same order $1/\sqrt{C}$ and the angle tends to $\pi/2$, hence we can choose the lift so that $\tilde{\go }\ra 0$.
\item As $C\ra 1^\pm $, the $T^{\rm Cas } $ tends to a finite limit, while $T^{\rm geod}\ra \infty $ and the angle $\ga $ tends to $0$ for most of the time along the closed orbit of the geodesic.
  Hence $\tilde{\go} \ra \infty $. .
  \item As $C \ra 0$, $T^{\rm Cas } $ tends to $\infty $ while $T^{\rm geod}$ is smooth. The angle $\ga $ is bounded below. Hence we can choose $\tilde{\go } \ra 0$. 
  \item As $C\ra -1$, $T^{\rm Cas } $ tends to $\sqrt{2}\pi  $ while $T^{\rm geod}$ is smooth. The angle $\ga $ tends to $0$, hence $\tilde{\go} \ra \infty $.
    
  \end{itemize}
  From all this, we see that the function $\tilde{\go} $ which is analytic is non constant on any interval, hence a dense set of values of $C$ for which we get periodic geodesics. 
\section{Appendices} 
\appendix
\section{Manifolds of fixed points}\label{app:inv}
Our goal is to prove the following
\begin{theo}\label{theo:inv}
  Let $Y_0 $ a compact submanifold of a manifold $X$ and $F_0:X\ra X $ a smooth map satisfying, $\forall y \in Y_0 $,
  $F_0(y)=y$ and $\ker (F'_0(y) -{\rm Id} )= T_y Y_0 $.
    and consider a smooth family of maps
  $F_\ge $. Then, for any $m$, there exists $\ge (m)>0$ so that,  for $|\ge| < \ge (m)$, there exists a $C^m$ manifold $Y_\ge $ depending smoothly  of $\ge $ globally invariant by $F_\ge $. 
 \end{theo}
In our paper this will be used with  $Y_0$ a circle embedded into a cylinder.

For the
proof, we use the  simplifying assumption that the normal bundle $TX_{|Y_0}/TY_0 $ is trivial. We can then reduce to the case where $X=(Y_0)_y\times \R^n_z $.

We will search $Y_\ge $ as a graph of a $C^m$ map $f_\ge $ from $Y_0$ into $\R^n$.
The invariance by $F_\ge(y,z)=(A_\ge (y,z), B_\ge (y,z) )$ writes
\[  B_\ge (y,f_\ge (y))=f_\ge (A_\ge (y,f_\ge (y)))    \]
Differentiating with respect to $\ge $ at $\ge =0$ gives
\[\forall y\in Y_0,  ({\rm Id}-B'_0 (y))\gd f(y)= \gd B (y,0)\]
where  the $\delta$'s are the derivatives w.r. to $\ge $ and  $B'_0 (y)$ is the derivative of $B_0$ with respect to $z$ at the point $(y,0)$ . Note that there is no derivatives of $A_\ge $ appearing because  $f_0=0$.  Hence the derivative of the righhandside w.r. to $\ge$
is $\gd f $. 
This can be solved with  $\gd f \in C^m (Y,\R^n )$ by the assumption of the Theorem.
Hence we can apply implicit function Theorem in the Banach space  $C^m (Y,\R^n )$  and conclude. 

\section{Poincaré-Birkhoff for twist maps}
The goal is to give a simple proof of the Poincaré-Birkhoff theorem for twist map of the annulus (see \cite{Go-01}). 
Let $A=(\R/\Z)_x \times [a,b]_y $ be an annulus equipped with some area form.
A smooth map $F=(X,Y):A\ra A $ preserving the boundaries of $A$ is called a {\it  twist map} if
\begin{enumerate}
  \item 
    $\pa X/\pa y >0 $
  \item $F$ is area preserving
    \item There
      exists  a choice $\tilde{F}=(\tilde{X},\tilde{Y})$ of a lift of $F $ to $\R\times [a,b]$ so that, for all $x\in \R$,  $\tilde{X}(x,a)<x $ and $\tilde{X}(x,b)>x $.
      \end{enumerate} 
Then
\begin{theo} \label{theo:PB}
  If $F$ is a twist map, it  admits a fixed point.
\end{theo}
The proof is as follows: for each $x\in \R/\Z$, the twist conditions (1)  and (3) implies that there exists a unique $y(x) $ so that $X(x,y(x))=x$. Moreover $x\ra y(x)$ is smooth.
Let us consider the curves which are the graphs of $x\ra y(x)$ and $x\ra  Y(x, y(x))$. The second curve is the image of the
first one by $F$. They have to cross because $F$ is area preserving (2), otherwise the image of the domain below the first curve will have an area smaller or larger than the area of that domain. 
Any intersection point is a fix point of $F$. 

\section{Perturbation of flows}\label{app:perturb}

Let us consider a vector field of the form
$\vec{V}=\vec{V_0}+\vec{R}$ and let $\phi_t^0 $ (resp. $\phi_t$) the  flows of $\overrightarrow{V_0}$ (resp. $\vec{V}$). 
We want to compare both flows. For that, we write $\phi_t (x)$ as $\phi_t (x)=\phi_t^0 (w(t,x) )$, with $w(0,x)=x$,  and will write a differential equation for
$w$.
We get the two equations for the integral curve $t\ra y(t)=\phi_t (x)$: 
\[  \frac{d}{dt} y(t)= \overrightarrow{V_0}(\phi_t^0 (w(t,x) ))+ \vec{R}(\phi_t^0 (w(t,x) )) \]
and
\[  \frac{d}{dt} y(t)= D\phi_t^0 (w(t,x) )\frac{d}{dt} w(t,x)+ \overrightarrow{V_0}(\phi_t^0 (w(t,x) )) \]
By identification of both equations, we get
\[ D\phi_t^0 (w(t,x) )\frac{d}{dt} w(t,x)=\vec{R}(\phi_t^0 (w(t,x) )) \]
and finally
\begin{equation}\label{equ:perturb}
  \frac{d}{dt} w(t,x)=  \left( D\phi_t^0 (w(t,x) )\right)^{-1}   \vec{R}(\phi_t^0 (w(t,x) )) \end{equation}
This imply that, if we have a weak control of $w $ and moreover we know $\phi_t^0$ and  the inverse of its differential, we get the closeness  of both flows. 

\section{Poincaré maps}  \label{app:poincare}

Let $\gg $ be a periodic orbit of a vector field $\vec{V}$  in a manifold  $M$. We choose a germ of hypersurface $S$ transverse to $\gg$ at some point $x_0$.
Then we can define a return map $F$ along $\gg$ which is a germ of map from $(S,x_0)$ into itself. 
That germ is independent of the choices of $x_0$ and $S$ up to conjugation by a germ of diffeomorphism.
The linearisation $L$ of $F$ at $x_0$ is hence also well defined up to conjugation.
The orbit $\gg$ is said to be non degenerate if $1$ is not an eigenvalue of $L$. In this case, we have the
\begin{prop} If $\gg$ is non degenerate and 
  if $\vec{V'}$ is close enough to $\vec{V}$, there is a close orbit $\gg'$ of $\vec{V'}$ which depends smoothly of $\vec{V'}$.
\end{prop}
The non degeneracy condition allows to  apply the inverse function theorem to the perturbation of the return map.

In case of an Hamiltonian vector field, the manifold $M $ to be considered is not the full phase space, but the energy shell containing $\gg$. The germ $S$ is then symplectic
and $F$ is a symplectic germ of diffeomorphism.

The Hamiltonian $g^\star $ is $g^\star =\xi^2+\eta ^2 $.
Let us show that $g^\star $ is ``integrable'' in the sense of Section \ref{sec:int}

\section{Geodesic flow on hyperbolic surfaces and Poincaré group of $M:=\Gamma \backslash {\rm PSL}_2 (\R)$}\label{app:hyp}
For this section, one can look at \cite{Bu-92}. 
We consider, for $\Gamma  \subset G={\rm PSL}_2(\R)$ a co-compact lattice with no elliptic elements, the compact smooth oriented hyperbolic surface  $N:=\Gamma \backslash \Hy $, where $G$ acts on $ \Hy $ by
$g(z)=(az+b )/(cz +d )$ with
\[ g = \left( \begin{matrix} a & b \\ c & d \end{matrix} \right) \]
The group $G$ acts also on the unit tangent bundle of $N$ obtained by taking the derivative of the previous action. The action is transitive with a trivial  isotropy group.
The map $g \ra Dg (v_0)$ where $v_0=(i, \pa_y)$ identifies $\Gamma \backslash G $ to the unit tangent bundle of $N$.
The $1$-parameter group ${\rm exp }(tZ )$ acting on the right on $G$ identifies then with the geodesic flow with speed $2$.
For each hyperbolic element $\gg \in \Gamma $, there exists an unique periodic geodesic $c $ of length $l_\gg=\log \gl $, whose lift to $\Hy $ is invariant by $\gg$.
We have  
$g {\rm exp }(\ha l_\gg Z) = \gamma g $. If $c $ is given, $\gg$ is determined up to conjugation in $\Gamma $. 

Given a group $H$, we denote by $C(H)$ the set of conjugacy classes of $H$. 
We have the following exact homotopy sequence associated to the fibration
$ \Gamma   \ra {\rm PSL}_2(\R) \ra \Gamma \backslash {\rm PSL}_2(\R) $: 
\[ \{ 1 \}\ra  \Z \ra \pi_1\left(\Gamma \backslash {\rm PSL}_2(\R) \right) \ra \Gamma  \ra \{ 1 \}  \]
The first arrow is an injective morphism whose image is the center of  $K:=\pi_1\left(\Gamma \backslash {\rm PSL}_2(\R) \right) $. The second arrow is surjective so that 
$K$ is an extension of $\Gamma $ by $\Z$. Concerning conjugacy classes  in $K$, there is an action of $\Z$ on $C(K)$  whose orbits are the fibers of the projection
$C(K)\ra C( \Gamma)$. 
In fact there is a canonical parametrization of $C(K)$ by $C(\Gamma )\times \Z$ as follows:
any element of $C(\Gamma )$ can be represented by a closed geodesic $\gg $ of $\Gamma \backslash \Hy $ which can be lifted as a perodic curve in $M$ by looking at the canonical lift
of $\gg $ to the unit tangent bundle. Then we look at the action of $\Z$ by composing this loop with a loop consisting in rotating the unit vector at a fixed point of $\Gamma \backslash \Hy$ of $2 \pi $. 
\bibliographystyle{plain}

\end{document}